\documentclass[letterpaper, 10 pt, conference]{ieeeconf}  

\IEEEoverridecommandlockouts                              

\usepackage{amsmath,amssymb,graphicx,fullpage,color,subcaption,wrapfig,epstopdf,diagbox, cuted, cite}
\usepackage[titlenumbered,ruled]{algorithm2e}
\usepackage[figuresright]{rotating}

\newcommand{\real}{\mathbb{R}}

\newcommand{\vo}[1]{\boldsymbol{#1}}

\newcommand{\Del}{\vo{\Delta}}

\title{\LARGE \bf
Robust LQR for Uncertain Discrete-Time Systems using Polynomial Chaos
}

\author{Vaishnav Tadiparthi$^{1}$ and Raktim Bhattacharya$^{2}$
\thanks{$^{1}$Vaishnav Tadiparthi is a graduate student at the Department of Aerospace Engineering Engineering, Texas A\&M University, College Station, TX 77843-3141. Email:
        {\tt\small vaishnavtv@tamu.edu}}%
\thanks{$^{2}$Raktim Bhattacharya is with the Faculty of Aerospace Engineering,  Texas A\&M University, College Station, TX 77843-3141. Email:
        {\tt\small raktim@tamu.edu}}%
}

\begin{document}

\maketitle
\thispagestyle{empty}
\pagestyle{empty}

\begin{abstract}
   In this paper, a polynomial chaos based framework for designing controllers for discrete time linear systems with probabilistic parameters is presented. Conditions for exponential-mean-square stability for such systems are derived and algorithms for synthesizing optimal quadratically stabilizing controllers are proposed in a convex optimization formulation. The solution presented is demonstrated on the derived discrete-time models of a nonlinear F-16 aircraft model trimmed at a set of chosen points.
\end{abstract}

\section{Introduction}

In this paper, we address the problem of designing linear quadratic regulators (LQRs) for discrete-time linear-time-invariant (LTI) systems with parametric uncertainty in the system matrices. Such systems can be defined as:
\begin{align}
   \vo{x}^{t+1}  = \vo{A(\Delta)}\vo{x}^t + \vo{B(\Delta)}\vo{u}^t
   \label{eqn:sysEq}
\end{align}
where the system matrices $\vo{A(\Delta)} \in \real^{n \times n}$ and $\vo{B(\Delta)} \in \real^{m \times n}$ are assumed to be affine functions in $\vo{\Delta} \in \mathcal{D}_{\Del} \subset \real^d$, representing the parameter space with a chosen probability density function $p(\Del)$ \cite{barmish1997uniform}.
The objective is to design a parameter-independent state-feedback law of the form $\vo{u} = \vo{Kx}$ that optimizes the closed-loop performance of the system in an LQR sense.

We solve the problem by \textit{maximizing the lower bound}\cite{willems1971least} on the cost-to-go using a reduced order model.
In \cite{daafouz2001parameter} and \cite{de1999new}, the authors described conditions for establishing robust stability in discrete time systems with time-varying parametric uncertainties.  Robust stabilization of similar systems with polytopic uncertainty by decoupling the Lyapunov and system matrices has also been addressed in the past \cite{zhang2007improved}.
However, the state of the art for handling systems like \eqref{eqn:sysEq} involve randomized algorithms \cite{tempo2012randomized}, a framework that requires large datasets to establish desired probabilistic guarantees \cite{kim2013wiener} and becomes computationally intractable for high dimensional parameter spaces.
Polynomial Chaos theory (PC) on the other hand, evokes deterministic algorithms that do not suffer from confidence issues, unlike randomized algorithms \cite{bhattacharya2014robust}.

We employ polynomial chaos theory as a means of numerical approximation \cite{xiu2002wiener} of the system and its uncertainties.
The PC framework is increasingly being used in the modeling of systems with random characteristics
 \cite{fisher2009linear}. In \cite{bhattacharya2019robust}, the authors presented a convex optimization formulation in the PC setting for the design of robust quadratic regulators for linear continuous systems \cite{stengel1991stochastic, polyak2001probabilistic} using deterministic algorithms.
Through this work, we expand on their contributions by extending one of the proposed approaches to discrete time systems with parametric uncertainties.
Furthermore, we derive an expression guaranteeing stability of such systems modeled using polynomial chaos theory.
With the conjunction of PC techniques and some useful results associated with Kronecker products \cite{bhattacharya2014robust}, we have truncated the dimension of the function space, thereby making the reduced order model more amenable to control design.

The paper is organized as follows: in section 2, we introduce the polynomial chaos framework that allows for the approximation techniques that follow. Using this framework, we present conditions for stability using standard Lyapunov theory. Formulation of the control objective and the solution methodology using linear matrix inequalities (LMIs) in a convex optimization problem is presented in section 3. The theorems are then applied to an illustrative control design problem obtained from a nonlinear F-16 longitudinal aircraft model \cite{stevens2015aircraft} and the paper concludes with a summary of the findings.
\section{Preliminaries}
\subsection{Polynomial Chaos Theory}
Polynomial Chaos (PC) theory is a deterministic framework to express the evolution of uncertainty in a dynamical system with probabilistic parameters.

 Using PC-based approximation techniques, $\vo{x}(t,\Del)$ can be expanded as:
\begin{align}
   \vo{x}^t(\Del) := \sum_{i=0}^{\infty}\vo{x}_i^t\phi_i(\Del)
\end{align}
where $\vo{x}_i^t$ are time-varying PC coefficients, and $\phi_i(\Del)$ are a known set of basis polynomials.
For the sake of brevity, we omit discussion over the nature of polynomials chosen, and encourage readers to refer to \cite{bhattacharya2019robust} for an in-depth commentary over the basis.
For computational tractability, we truncate the PC expansion to a finite number of terms, i.e.,
\begin{align}
   \vo{x}^t(\Del) \approx \vo{\hat{x}}^t(\Del) := \sum_{i=0}^{N}\vo{x}_i^t\phi_i(\Del)
   \label{eqn:xApproxSum}
\end{align}

\subsection{Surrogate System Modeling with Polynomial Chaos}
Define $\vo{\Phi}$ to be:
\begin{align}
   \vo{\Phi}(\vo{\Delta}) :=& (\phi_0(\vo{\Delta}), \phi_1(\vo{\Delta}) {}, \hdots, {} \phi_N(\vo{\Delta})), \textrm{and} \\
   \vo{\Phi}_n(\vo{\Delta}) :=& \vo{\Phi}(\vo{\Delta}) \otimes \vo{I}_n,
\end{align}

where $\vo{I}_n \in \real^{n \times n}$ is the identity matrix. Now, define matrix $\vo{X} \in \real^{n \times (N+1)}$, with PC coefficients $\vo{x}_i \in \real^n$ as
\begin{align}
   \vo{X} = \vo{[x}_0, \hdots, \vo{x}_N].
\end{align}
Therefore, \eqref{eqn:xApproxSum} can be written compactly:
\begin{align}
   \vo{\hat{x}}^t(\vo{\Delta}) := \vo{X}^t\vo{\Phi}(\vo{\Delta}).
\end{align}
Further vectorizing,
\begin{align}
   \vo{\hat{x}}^t(\vo{\Delta}) &\equiv \textrm{vec}(\vo{\hat{x}}^t(\vo{\Delta})) \nonumber\\
   &= \textrm{vec}(\vo{X}^t\vo{\Phi}(\vo{\Delta})) \nonumber\\
   &= \textrm{vec}(\vo{I}_n \vo{X}^t \vo{\Phi}(\vo{\Delta})) \nonumber\\
   &= \vo{\Phi}_n^T(\vo{\Delta})\vo{x}_{pc}^t
\end{align}
where $\vo{x}_{pc}^t := \textrm{vec}(\vo{X}^t)$
The error due to approximation can be expressed as:
\begin{align}
   \vo{e}^t(\Del) &:= \vo{\hat{x}}^{t+1} - \vo{A}(\Del)\vo{\hat{x}}^t - \vo{B}(\Del)\vo{K}\vo{\hat{x}}^t \nonumber\\
   &= \vo{\Phi}_n^T(\vo{\Delta})\vo{x}_{pc}^{t+1} - \vo{A}(\Del)\vo{\Phi}_n^T(\vo{\Delta})\vo{x}_{pc}^t - \nonumber\\ & \quad \vo{B}(\Del)\vo{K}\vo{\Phi}_n^T(\vo{\Delta})\vo{x}_{pc}^t
\end{align}

To minimize the error in the $\mathcal{L}_2$ sense, we need the projection of the expected value of this approximation error on
the basis function to be zero, i.e.,
\begin{align}
   \mathbb{E}[\vo{e}^t(\Del)\phi_i(\Del)] = 0 \textrm{ for } i = 0,1,\hdots,N
\end{align}

Upon simplification, we arrive at $n(N+1)$ \textit{deterministic} ordinary differential equations in $\vo{x}_{pc}^t$.
\begin{align}
   \vo{x}_{pc}^{t+1} = (\mathbb{E}[\vo{\Phi}(\Del)  \vo{\Phi}^T(\Del)] \otimes \vo{I}_n)^{-1} \times \nonumber\\
   (\mathbb{E}[(\vo{\Phi}(\Del)  \vo{\Phi}^T(\Del)) \otimes \vo{A}(\Del)] + \nonumber\\
   \mathbb{E}[(\vo{\Phi}(\Del)  \vo{\Phi}^T(\Del)) \otimes (\vo{B}(\Del)\vo{K})])\vo{x}_{pc}^t
	\label{eqn:pcODE}
\end{align}

\section{Controller Synthesis}
The design objective is to determine an optimal state-feedback law of the form $\vo{u} = \vo{Kx}$, to optimize the closed-loop performance in an LQR sense.

Therefore, the closed-loop system with control $\vo{u} = \vo{Kx}$ is given by,
\begin{align}
   \vo{x}^{t+1} = (\vo{A(\Delta) + B(\Delta)K})\vo{x}^t
   \label{eqn:sysEqCL}
\end{align}

The system in \eqref{eqn:sysEqCL} is infinite-dimensional with respect to $\vo{\Delta}$. We present a formulation that reduces the infinite-dimensional system using PC expansion and develop an optimization problem with the reduced-order system. The approximate problem is then solved exactly.

\subsection{Stability}
An almost-sure (a.s.) stability analysis of linear systems from a polynomial chaos framework is prohibitive, except for some special cases \cite{dutta2010nonlinear, ghanem2003stochastic}. However, the equivalence of a.s. stability to exponential-mean-square (EMS) stability for LTI systems \cite{chen1995linear} means that we can use EMS  analysis, a far more favorable framework based on moments, to examine our system in the PC setting.

\textit{Derivation}:
Define a parameter-dependent Lyapunov function $V(\vo{x}) := \vo{x}^T \vo{P}(\Del) \vo{x}$.\\
Represent $\vo{P(\Delta)}$ using homogeneous polynomial-based parameter-dependent quadratic functions, i.e.,
\begin{align}
	\vo{P(\Delta)} := {\vo{\Phi_n^T (\Delta)}} \bar{\vo{P}} \vo{\Phi_n (\Delta)}, \, \bar{\vo{P}} \in \mathcal{S}_{++}^{n(N+1)}
\end{align}
If $\vo{P}(\Del) > \vo{0},$ $\forall\Del $ $\iff \bar{\vo{P}} > \vo{0}$. Partitioning $\vo{P}$ as:
\begin{align*}
  \bar{\vo{P}} := \begin{bmatrix}
    \vo{P}_{00} & \hdots  & \vo{P}_{0N} \\
    \vdots & {} & \vdots \\
    \vo{P}_{N0} & \hdots  & \vo{P}_{NN}
\end{bmatrix}
\end{align*}
where $\vo{P}_{ij} = \vo{P}_{ji} \in \real^{n \times n}$. Therefore,
\begin{align}
  \vo{P(\Delta)} &:= {\vo{\Phi_n^T (\Delta)}} \bar{\vo{P}} \vo{\Phi}_n (\Del) \nonumber \\
  &= \begin{bmatrix}
    \phi_0(\Del)\vo{I}_n & \hdots & \phi_N  (\Del)  \vo{I}_n
\end{bmatrix} \times \nonumber \\
& \begin{bmatrix}
  \vo{P}_{00} & \hdots  & \vo{P}_{0N} \\
  \vdots & {} & \vdots \\
  \vo{P}_{N0} & \hdots  & \vo{P}_{NN}
\end{bmatrix} \times \begin{bmatrix}
  \phi_0(\Del)\vo{I}_n \\ \vdots \\ \phi_N  (\Del)  \vo{I}_n
\end{bmatrix} \nonumber \\
&= \sum_{ij} \phi_i (\Del) \phi_j (\Del) \vo{P}_{ij}
\end{align}
Noting that $\vo{P}(\Del)$ is symmetric,
\begin{align*}
  \vo{P}_{ij} = \vo{P}_{ji}.
\end{align*}

\textbf{Theorem}:
The system given in \eqref{eqn:sysEq} is EMS-stable if $\exists$ $\vo{\bar{P}} = \vo{\bar{P}}^T > 0$ such that $\mathbb{E} [\vo{W}_{pc}] < 0$ where $\mathbb{E} [\vo{W}_{pc}]$ is given element-wise by the expression below.
\begin{align}
        \mathbb{E} [\vo{W}_{pc_{il}}] = &\sum_{j,k = 0}^{N} \sum_{m,n = 0}^{5} \left[ \mathbb{E}[\phi_i \phi_j \phi_k \phi_m \phi_n \phi_l] \times \right. \nonumber\\
         & \left. \vo{A}_m^T \vo{P}_{jk} \vo{A}_n \right] 
        - \sum_{j,k = 0}^{N} \left[ \mathbb{E}[\phi_i \phi_j \phi_k \phi_l] \vo{P}_{jk} \right] \nonumber\\
        &\textrm{for } i, l=0,1,\hdots,N
  \label{eqn:stab}
\end{align}

\textbf{\textit{Proof}}: Examining the expression $\mathbb{E}[\Delta V]$,
\begin{align*}
  \mathbb{E}[\Delta V] &= V^{t+1} - V^{t} \\
  &= (\vo{x}^{t+1})^T \vo{P}(\Del) \vo{x}^{t+1} - \vo{x}^{tT} \vo{P}(\Del) \vo{x}^t \\
  &= (\vo{A}(\Del)  \vo{x}^t)^T \vo{P}(\Del)(\vo{A}(\Del) \vo{x}^t) -  \vo{x}^{tT} \vo{P}(\Del) \vo{x}^t \\
  &=  \vo{x}^{tT}( \vo{A}^T (\Del)\vo{P}(\Del) \vo{A} (\Del)- \vo{P}(\Del)) \vo{x}^t
\end{align*}
Dropping the $t$ superscript for notational convenience,
\begin{align*}
  \mathbb{E}[\Delta V] = \vo{x}^T (\vo{A}^T(\Del) \vo{P}(\Del) \vo{A} (\Del)- \vo{P}(\Del))\vo{x}
\end{align*}
For stability,
  $\mathbb{E}[\Delta V] < 0$.
Substituting $\vo{x}^t = \vo{\Phi}_n^T(\vo{\Delta})\vo{x}_{pc}^t$,
\begin{align}
  \mathbb{E}[\Delta V] = \vo{x}_{pc}^T \mathbb{E}[\vo{\Phi}_n(\vo{\Delta})(\vo{A}^T (\Del)\vo{P}(\Del) \vo{A}(\Del) - \nonumber\\ \vo{P}(\Del)) \vo{\Phi}_n^T (\Del)] \vo{x}_{pc}
\end{align}
Let:
\begin{align*}
        \vo{W}_{pc} = \vo{\Phi}_n(\vo{\Delta})(\vo{A}^T (\Del)\vo{P}(\Del) \vo{A}(\Del) - \vo{P}(\Del)) \vo{\Phi}_n^T (\Del)
\end{align*}
Expanding the expressions:
\begin{align*}
        \vo{\Phi}_n (\Del) = (\vo{\Phi} (\Del) \otimes \vo{I}_n), \quad
        \vo{A}(\Del) = \sum_{m=0}^{\textrm{nOrd}} \vo{A}_m \phi_m(\Del)
\end{align*}
'$\textrm{nOrd}$' is the order of polynomials chosen to approximate the uncertain system matrices $\vo{A}$ and $\vo{B}$.
Using the standard properties of Kronecker products, and switching to the index notation, we obtain the equation \eqref{eqn:stab}. Therefore, this reduces to:
\begin{align*}
        \mathbb{E} [\vo{W}_{pc_{il}}] = &\sum_{j,k = 0}^{N} \sum_{m,n = 0}^{5} \left[ \mathbb{E}[\phi_i \phi_j \phi_k \phi_m \phi_n \phi_l] \times \right. \nonumber\\
        & \left. \vo{A}_m^T \vo{P}_{jk} \vo{A}_n \right] 
        - \sum_{j,k = 0}^{N} \left[ \mathbb{E}[\phi_i \phi_j \phi_k \phi_l] \vo{P}_{jk} \right] \nonumber\\
        &\textrm{for } i, l=0,1,\hdots,N
\end{align*}
Note that for each $i$ and $l$, $\mathbb{E}[V_{pc_{il}}] \in \real^{n \times n}$.

\subsection{LQR}
Assuming the system given by \eqref{eqn:sysEq} is EMS-stable, we consider the synthesis of an optimal state-feedback control gain that minimizes a quadratic cost, i.e., a fixed parameter-independent gain $\vo{K}$ that minimizes
\begin{align}
	\mathbb{E} \left[ \sum_{t=0}^{\infty} (\vo{x}^t)^T \vo{Q} \vo{x}^t + (\vo{u}^t)^T \vo{R} \vo{u}^t)\right]
\end{align}
subject to $\vo{u}^t = \vo{K}\vo{x}^t$ and dynamics given by \eqref{eqn:sysEqCL}.


For the dynamical system given in \eqref{eqn:sysEqCL} with controller $\vo{u} = \vo{Kx}$, the solution is:
\begin{align}
	\vo{x}^t(\Del) = (\vo{A_c})^{t}(\Del)\vo{x}_0
\end{align}
where $\vo{A}_c(\Del) := \vo{A}(\Del) + \vo{B}(\Del)\vo{K}$. The cost-to-go from initial condition $x_0$ is therefore,
\begin{gather}
	\vo{J}(\vo{x}_0) = \mathbb{E}\left[ \sum_{t=0}^{\infty} (\vo{x}^t)^T \vo{Q} \vo{x}^t + (\vo{u}^t)^T \vo{R} \vo{u}^t)\right] \nonumber\\
							= \vo{x_0}^T \mathbb{E} \underbrace{\left[  \sum_{t=0}^{\infty} (\vo{A}_c^{t}(\Delta))^T(\vo{Q} + \vo{K}^T\vo{RK}) \vo{A}_c^{t}(\Del) \right]}_{\vo{:=P(\Delta)}} \vo{x_0}
\end{gather}
for some $\vo{P(\Delta)} : \real^d \mapsto \mathcal{S}_{++}^n$. Therefore, the cost-to-go is a quadratic fucntion of states.

\subsection{Model Reduction}
Rewriting \eqref{eqn:pcODE} using properties associated with Kronecker products \cite{bhattacharya2014robust}, we obtain the reduced order approximation:
\begin{align}
   \vo{x}_{pc}^{t+1} = (\vo{A}_{pc} + \vo{B}_{pc}\vo{\mathcal{K}})\vo{x}_{pc}^{t}
\end{align}
where $\vo{\mathcal{K}} = \vo{I}_{N+1} \otimes \vo{K}$, and
\begin{align}
   \vo{A}_{pc} &:= \vo{\Phi_K} 
	\mathbb{E}[(\vo{\Phi}(\Del)  \vo{\Phi}^T(\Del)) \otimes \vo{A}(\Del)] \\
   \vo{B}_{pc} &:= \vo{\Phi_K}  
	\mathbb{E}[(\vo{\Phi}(\Del)  \vo{\Phi}^T(\Del)) \otimes \vo{B}(\Del)]
\end{align}
where $\vo{\Phi_K} := (\mathbb{E}[\vo{\Phi}(\Del)  \vo{\Phi}^T(\Del)] \otimes \vo{I}_n)^{-1}$.

The modified cost function is:
\begin{align}
   \vo{J} := \sum_{t=0}^{\infty}(\vo{x}_{pc}^t)^T (\vo{Q}_{pc} + \vo{\mathcal{K}}^T \vo{R}_{pc} \vo{\mathcal{K}})\vo{x}_{pc}^t
\end{align}
where:
\begin{align}
   \vo{Q}_{pc} &:= \mathbb{E}[\vo{\Phi}_n(\Del) \vo{Q} \vo{\Phi}_n^T(\Del)] \\
   \vo{R}_{pc} &:= \mathbb{E}[\vo{\Phi}_n(\Del) \vo{R} \vo{\Phi}_n^T(\Del)]
\end{align}
Defining the Lyapunov function in $\vo{x}_{pc}^t$, i.e.,
\begin{align}
   V(\vo{x}_{pc}^t) := (\vo{x}_{pc}^t)^T \vo{P}_{pc} \vo{x}_{pc}^t
\end{align}
where $\vo{P}_{pc} \in \mathcal{S}^{n(N+1)}_{++}$. Dropping the $t$ superscript for notational convenience, we get:
\begin{align}
   &[\Delta V] +  [\vo{x}_{pc}^T (\vo{Q}_{pc} + \vo{\mathcal{K}}^T \vo{R}_{pc} \vo{\mathcal{K}}) \vo{x}_{pc} ]  \nonumber\\
      &= \vo{x}_{pc}^T [(\vo{A}_{pc} + \vo{B}_{pc} \vo{\mathcal{K}})^T \vo{P}_{pc}  (\vo{A}_{pc} + \vo{B}_{pc}\vo{\mathcal{K}}) \nonumber\\
      & - \vo{P}_{pc} + \vo{Q}_{pc} + \vo{\mathcal{K}}^T \vo{R}_{pc} \vo{\mathcal{K}}] \vo{x}_{pc}
      \nonumber \\
    & = \vo{x}_{pc}^T [ \vo{A}_{pc}^T \vo{P}_{pc} \vo{A}_{pc} - \vo{P}_{pc} + \vo{Q}_{pc}\nonumber\\
    & - \vo{A}_{pc}^T \vo{P}_{pc} \vo{B}_{pc} (\vo{R}_{pc} + \vo{B}_{pc}^T \vo{P}_{pc} \vo{B}_{pc})^{-1} \vo{B}_{pc}^{T} \vo{P}_{pc} \vo{A}_{pc} \nonumber\\
     &+ (\vo{\mathcal{K}} + (\vo{R}_{pc} + \vo{B}_{pc}^T \vo{P}_{pc} \vo{B}_{pc})^{-1} \vo{B}_{pc}^T \vo{P}_{pc} \vo{A}_{pc})^T \times \nonumber\\
      &(\vo{R}_{pc} + \vo{B}_{pc}^T \vo{P}_{pc} \vo{B}_{pc}) \times \nonumber \\
     & (\vo{\mathcal{K}} + (\vo{R}_{pc} + \vo{B}_{pc}^T \vo{P}_{pc} \vo{B}_{pc})^{-1} \vo{B}_{pc}^T \vo{P}_{pc} \vo{A}_{pc})] \vo{x}_{pc}
      \label{eqn:redOrderEqn}
\end{align}

We cannot determine $\vo{K}$ by setting $\vo{\mathcal{K}} := \vo{R}_{pc} + \vo{B}_{pc}^T \vo{P}_{pc} \vo{B}_{pc})^{-1} \vo{B}_{pc}^T \vo{P}_{pc} \vo{A}_{pc}$, unless we are assured that all of $\vo{A}_{pc}$, $\vo{B}_{pc}$, and $\vo{P}_{pc}$
have block diagonal structures with repeating blocks, thus leading to our formulation:

\textbf{Theorem:}
 The optimal controller $\vo{K}$ and $\vo{P}_{pc}$ are obtained as solutions to the following optimization problems:
\begin{gather}
	\underset{\vo{P}_{pc} \in {\mathcal{S}}_{++}^{n(N+1)}}{\textrm{max}} \textrm{tr } \vo{P}_{pc}, \nonumber
        \\  \textrm{subject to} \nonumber\\
	\begin{bmatrix}
		\vo{A}_{pc}^T \vo{P}_{pc} \vo{A}_{pc} + \vo{Q}_{pc} - \vo{P}_{pc} & \vo{A}_{pc}^T \vo{P}_{pc} \vo{B}_{pc} \\
		\vo{B}_{pc}^T \vo{P}_{pc} \vo{A}_{pc} & \vo{R}_{pc} + \vo{B}_{pc}^T \vo{P}_{pc} \vo{B}_{pc}
	\end{bmatrix} \nonumber\\
	\geq \vo{0}
        \label{eqn:thmP} \\
\textrm{and} \nonumber\\
	\underset{\vo{X} \in \vo{\mathcal{S}}_{++}^{n(N+1)}, \vo{K} \in \real^{m \times n}}{\textrm{min}} \textrm{tr } \vo{X}, \nonumber
        \\  \textrm{subject to} \nonumber\\
   \begin{bmatrix}
      \vo{H}_1 & \vo{H}_2 \\
      \vo{H}_2^T & \vo{H}_4
   \end{bmatrix}
	\geq \vo{0}
        \label{eqn:thmK}
\end{gather}
where
\begin{flalign*}
   \vo{H}_1 &= \vo{X} \\
   \vo{H}_2 &= (\vo{\mathcal{K}} + (\vo{R}_{pc} + \vo{B}_{pc}^T \vo{P}_{pc} \vo{A}_{pc})^{-1} \vo{B}_{pc}^T \vo{P}_{pc} \vo{A}_{pc})^T \\
   \vo{H}_4 &= (\vo{R}_{pc} + \vo{B}_{pc}^T \vo{P}_{pc} \vo{B}_{pc})^{-1}
\end{flalign*}


\textit{\textbf{Proof:}}
Following the approach similar to the one described in \cite{bhattacharya2019robust} wherein they obtained the optimal solution by maximizing the lower bound, i.e., setting \eqref{eqn:redOrderEqn} $\geq 0$.
\begin{align}
	&\vo{A}_{pc}^T \vo{P}_{pc} \vo{A}_{pc} \nonumber\\
   & - \vo{A}_{pc}^T \vo{P}_{pc} \vo{B}_{pc} (\vo{R}_{pc} + \vo{B}_{pc}^T \vo{P}_{pc} \vo{B}_{pc})^{-1} \vo{B}_{pc}^T \vo{P}_{pc} \vo{A}_{pc} \nonumber\\
   & - \vo{P}_{pc} + \vo{Q}_{pc} \nonumber\\
	& + (\vo{\mathcal{K}} + (\vo{R}_{pc} + \vo{B}_{pc}^T \vo{P}_{pc} \vo{A}_{pc})^{-1} \vo{B}_{pc}^T \vo{P}_{pc} \vo{A}_{pc})^T \times \nonumber\\
   &(\vo{R}_{pc} + \vo{B}_{pc}^T \vo{P}_{pc} \vo{B}_{pc} ) \times \nonumber \\
   &(\vo{\mathcal{K}} + (\vo{R}_{pc} + \vo{B}_{pc}^T \vo{P}_{pc} \vo{A}_{pc})^{-1} \vo{B}_{pc}^T \vo{P}_{pc} \vo{A}_{pc}) \geq \vo{0} \label{eqn:maxLB}\\
	&\implies \Delta V \geq -\vo{x}_{pc}^T (\vo{Q}_{pc} + \vo{\mathcal{K}}^T \vo{R}_{pc} \vo{\mathcal{K}}) \vo{x}_{pc} \nonumber
\end{align}
Summing over $[0, \infty]$, we obtain:
\begin{align}
   \sum_{t=0}^{\infty} \Delta V^t \geq  -  \sum_{t=0}^{\infty} (\vo{x}_{pc}^t)^T (\vo{Q}_{pc} + \vo{\mathcal{K}}^T \vo{R}_{pc} \vo{\mathcal{K}})(\vo{x}_{pc}^t)
\end{align}
which is equivalent to:
\begin{align}
   \lim_{t \to \infty}V(\vo{x}_{pc}^t) - V(\vo{x}_{pc}^0) \geq \nonumber\\
	  - \sum_0^{\infty} \vo{x}_{pc}^T (\vo{Q}_{pc} + \vo{\mathcal{K}}^T \vo{R}_{pc} \vo{\mathcal{K}}) \vo{x}_{pc}
\end{align}
Since the closed loop system is EMS stabilizable,
\begin{align}
   \lim_{t \to \infty}V(\vo{x}_{pc}^t) = 0,
\end{align}
thus implying:
\begin{align}
   V(\vo{x}_{pc}^0) \leq  -  \sum_{t=0}^{\infty} \vo{x}_{pc}^T (\vo{Q}_{pc} + \vo{\mathcal{K}}^T \vo{R}_{pc} \vo{\mathcal{K}}) \vo{x}_{pc}
\end{align}

The inequality \eqref{eqn:maxLB} is nonconvex and is relaxed as:
\begin{align}
   &(\vo{\mathcal{K}} + (\vo{R}_{pc} + \vo{B}_{pc}^T \vo{P}_{pc} \vo{A}_{pc})^{-1} \vo{B}_{pc}^T \vo{P}_{pc} \vo{A}_{pc})^T  \times \nonumber\\
   &(\vo{R}_{pc} + \vo{B}_{pc}^T \vo{P}_{pc} \vo{B}_{pc})  \times \nonumber\\
   &(\vo{\mathcal{K}} + (\vo{R}_{pc} + \vo{B}_{pc}^T \vo{P}_{pc} \vo{A}_{pc})^{-1} \vo{B}_{pc}^T \vo{P}_{pc} \vo{A}_{pc}) \geq \vo{0}
\end{align}
for any $\vo{K}$ (or $\vo{\mathcal{K}}$). Therefore, the inequality transforms:
\begin{align}
   & - \vo{A}_{pc}^T \vo{P}_{pc} \vo{B}_{pc} (\vo{R}_{pc} + \vo{B}_{pc}^T \vo{P}_{pc} \vo{B}_{pc})^{-1} \ \vo{B}_{pc}^T \vo{P}_{pc} \vo{A}_{pc} \nonumber\\
   & + \vo{A}_{pc}^T \vo{P}_{pc} \vo{A}_{pc} - \vo{P}_{pc} + \vo{Q}_{pc} \geq \vo{0}
   \label{eqn:thmPpc}
\end{align}
Subsequently, determine a $\vo{K}$ that minimizes
\begin{align}
        (\vo{\mathcal{K}} + (\vo{R}_{pc} + \vo{B}_{pc}^T \vo{P}_{pc} \vo{A}_{pc})^{-1} \vo{B}_{pc}^T \vo{P}_{pc} \vo{A}_{pc})^T \times \nonumber\\
        (\vo{R}_{pc} +  \vo{B}_{pc}^T \vo{P}_{pc} \vo{B}_{pc}) \times \nonumber\\
        (\vo{\mathcal{K}} + (\vo{R}_{pc} + \vo{B}_{pc}^T \vo{P}_{pc} \vo{A}_{pc})^{-1} \vo{B}_{pc}^T \vo{P}_{pc} \vo{A}_{pc})
        \label{eqn:thmKpc}
\end{align}
Equations \eqref{eqn:thmPpc} and \eqref{eqn:thmKpc} can be expressed as the LMIs in \eqref{eqn:thmP} and \eqref{eqn:thmK}.


\section{Example and Results}
The plant considered here is a longitudinal F-16 aircraft model. The states of the system are $V$(ft/s), angle-of-attack $\alpha$(rad), and pitch rate $q$(rad/s). The control variables are thrust $T$ (lb) and elevator angle $\delta_e$(deg). The nonlinear model is trimmed at velocities from 400 ft/s to 900 ft/s in increments of 100 ft/s, and at an altitude of 10000 ft.
The objective is to design a fixed gain $\vo{K}$ that is able to regulate perturbations in the linear plant about various equilibrium points.

The vehicle is trimmed using a constraned nonlinear least-squares optimization, solved using sequential quadratic programming with \textit{fmincon} in MATLAB. The following constraints are imposed on the state and control magnitudes to account for aerodynamic and actuator limits:

\begin{center}
  \begin{tabular}{l r}
    Trim at $V_{trim}$: & $V_{trim} \leq V \leq V_{trim}$ \\
    Validity of aerodynamic data: & $-20^{\circ} \leq \alpha \leq 40^{\circ}$ \\
    Steady-level flight: & $0 \leq \theta - \alpha \leq 0$\\
    Steady-level flight: & $0 \leq q \leq 0$\\
    Thrust limits: & $1000 \leq T \leq 19000$ \\
    Elevator limits: & $-25^{\circ} \leq \delta_e \leq 25^{\circ}$\\
  \end{tabular}
\end{center}

\begin{figure}[!ht]
\begin{subfigure}[t]{0.23\textwidth}
   \includegraphics[width=\textwidth]{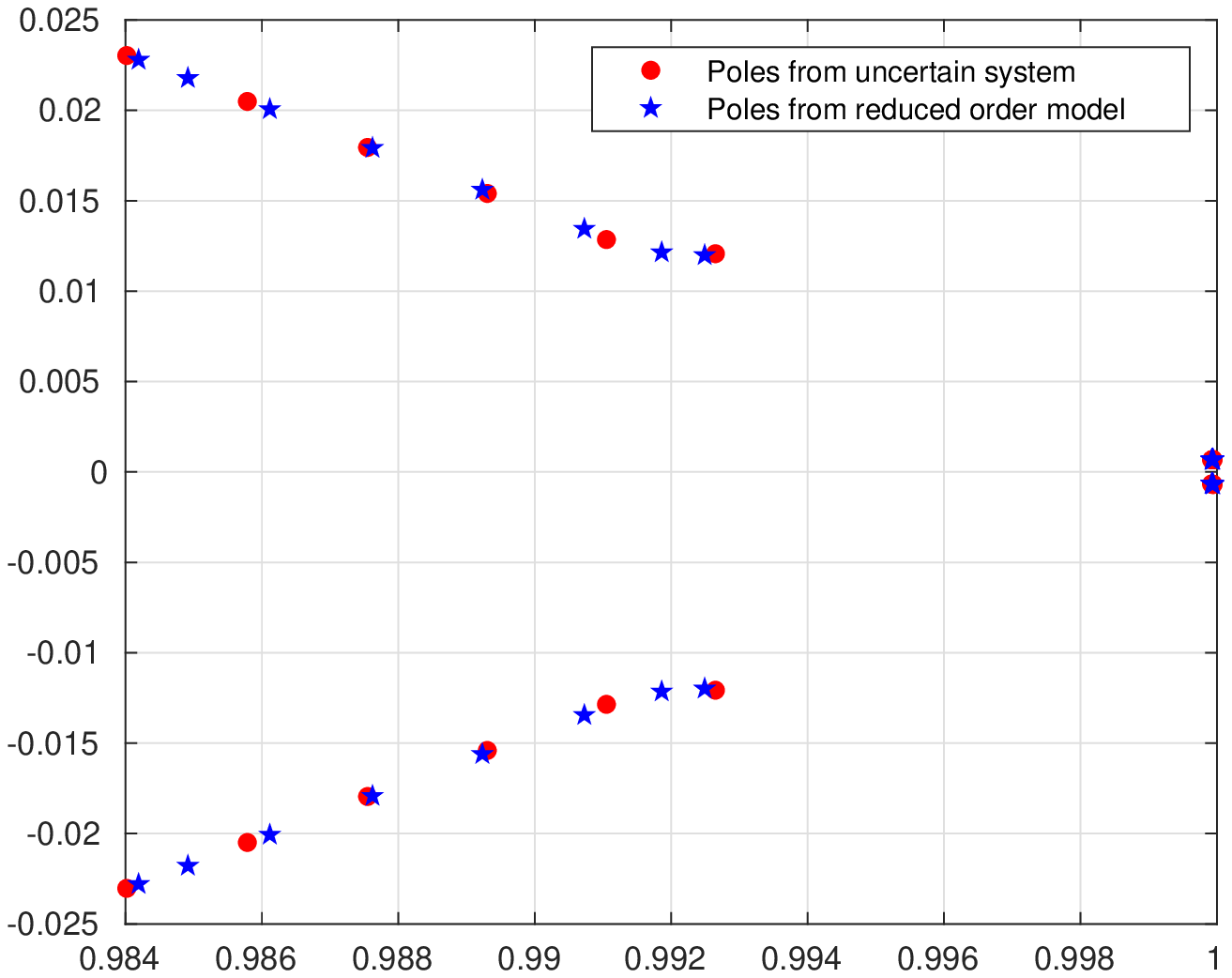}
   \caption{Open-loop poles of the uncertain model and the reduced-order model with 7th-order approximation.}
   \label{olPoles}
   \end{subfigure}\hfill
   \begin{subfigure}[t]{0.23\textwidth}
    \includegraphics[width=\textwidth]{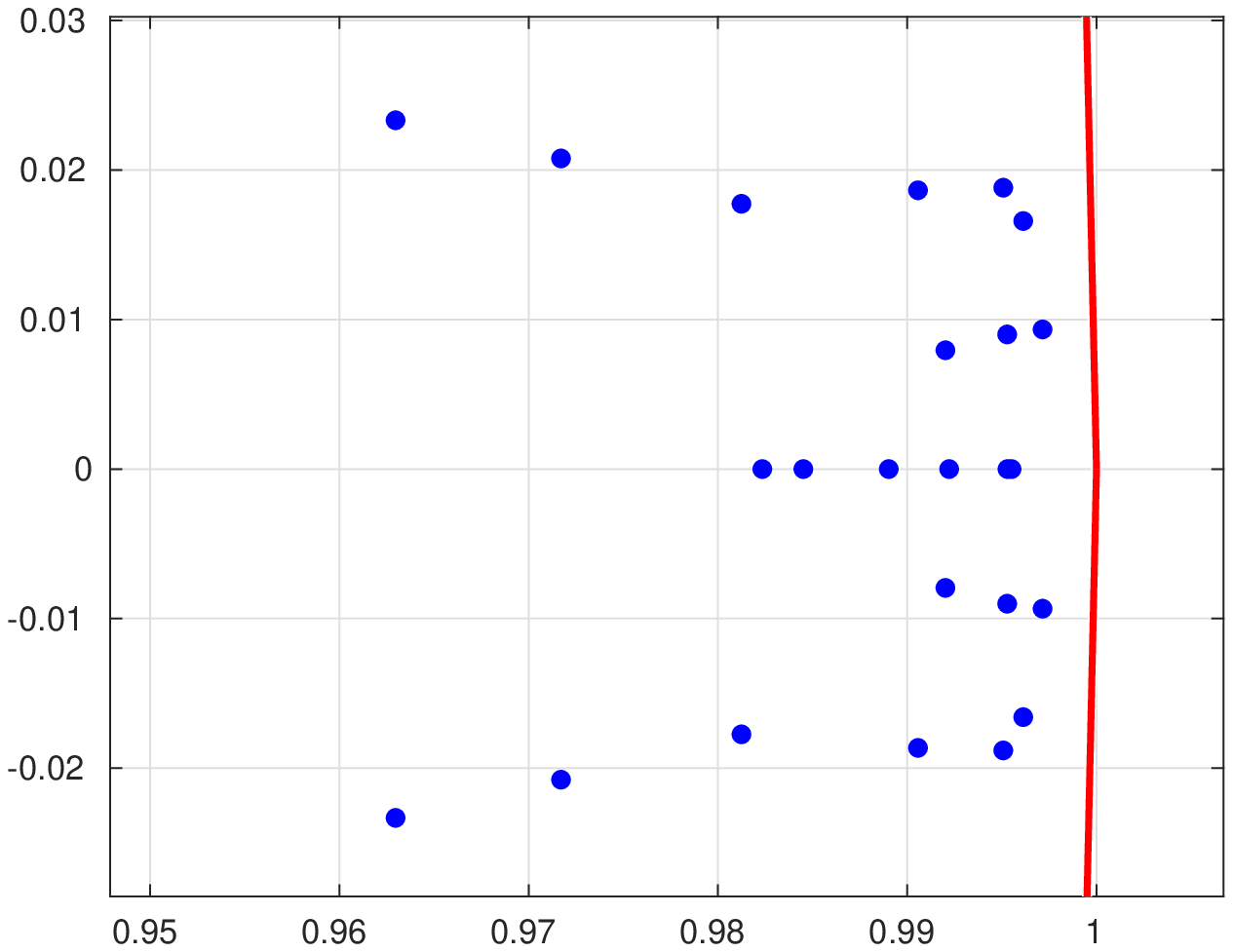}
   \caption{Closed-loop poles for the reduced order model.}
   \label{clLoopPoles}
   \end{subfigure}
   \caption{Open and closed-loop poles.}
\end{figure}

Linear models at each of these flight velocities were obtained using MATLAB's \textit{linmod} command. Subsequently, the discrete-time linear models were derived using the \texttt{c2d} command. Figure \ref{olPoles} shows the open loop poles of the discretized linear system. At $V_{trim} = 400$ ft/s, the system is marginally stable. At higher velocities, the poles are all comfortably inside the unit circle.

\begin{figure}[!ht]
\begin{subfigure}[t]{0.23\textwidth}
   \includegraphics[width=\textwidth]{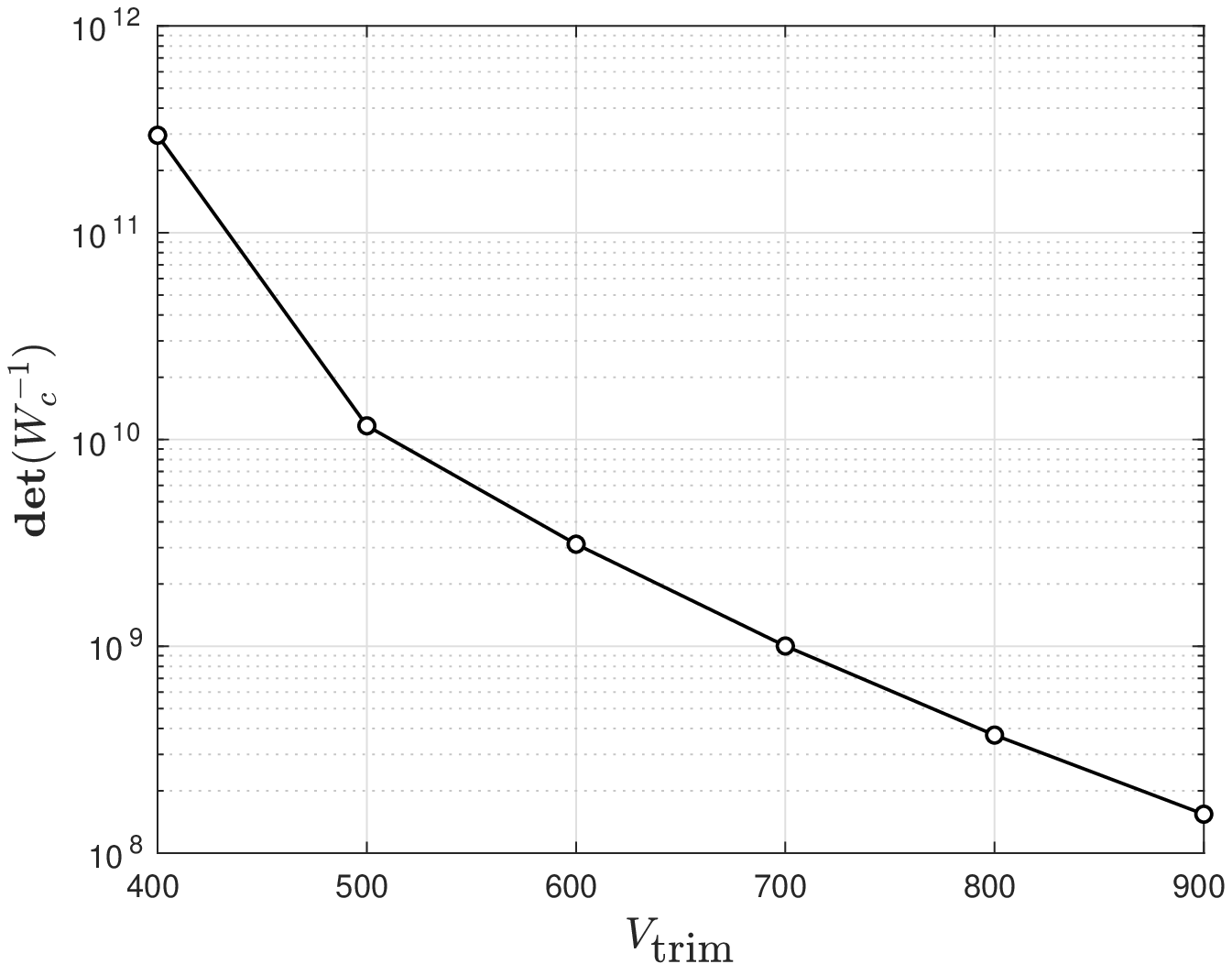}
   \caption{Degree of controllability of the system as a function of $\Del$. Note that the system is less controllable for smaller values of $\Del$.}
   \label{ctrbEnergy}
	\end{subfigure}\hfill
   \begin{subfigure}[t]{0.23\textwidth}
   \includegraphics[width=\textwidth]{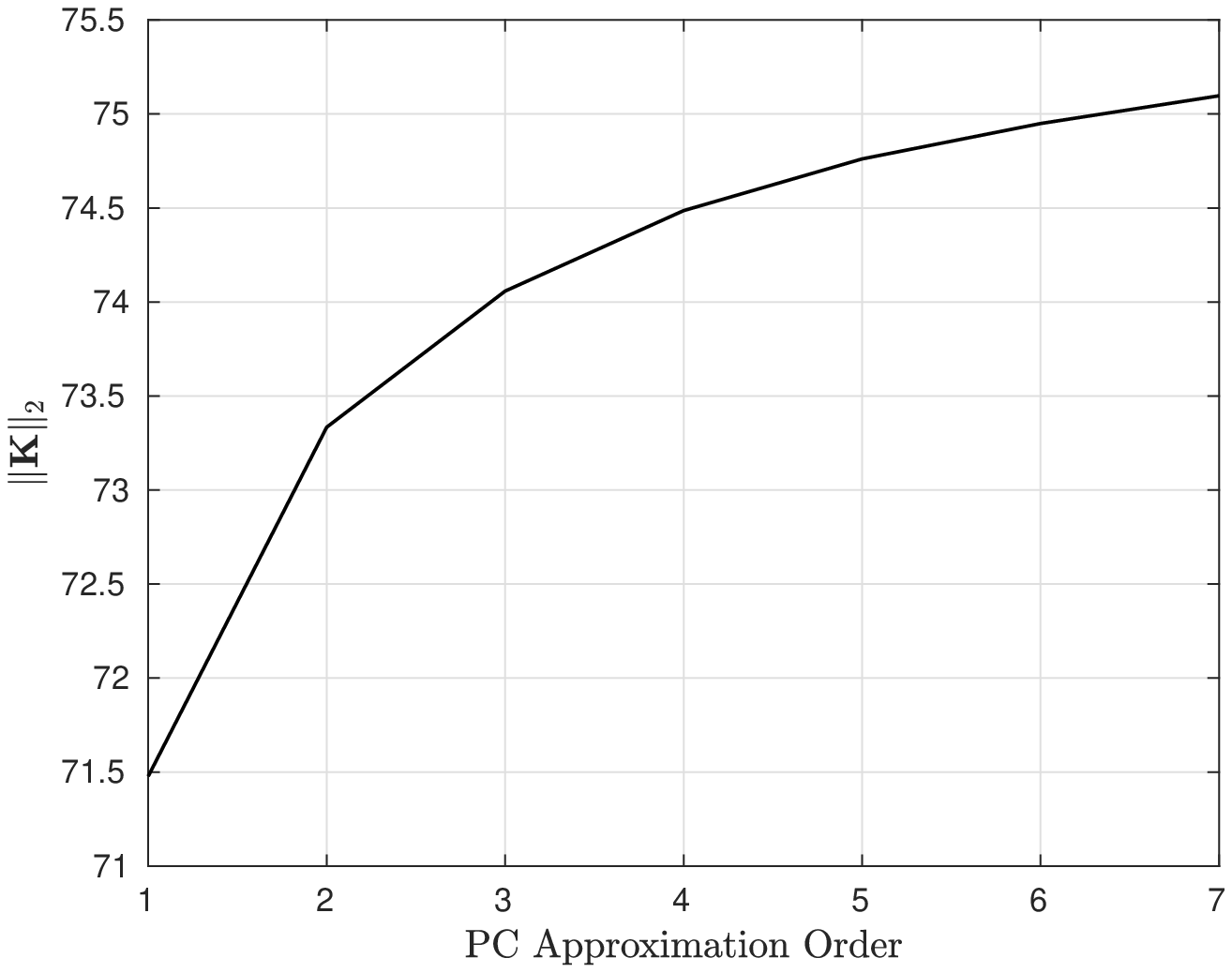}
   \caption{Variation of $||\vo{K}||_2$ with PC approximation order.}
   \label{kGain}
   \end{subfigure}
   \caption{Controllability and controller gain.}
\end{figure}

Figure \ref{ctrbEnergy} shows $\textbf{det}(W_c^{-1})$ where $W_c$ is the controllability Gramian, for various values of $V_{trim}$ for which the open-loop system is stable. This Gramian helps gauge the energy required to move the system around the state space. The determinant is a scalar measure of the Gramian, and figure \ref{ctrbEnergy}  shows that the energy increases for lower values of $V_{trim}$, indicating that closed-loop performance at these velocities degrades.

We scale $V$ and represent it as $\Del$, where
\begin{align*}
  \Del = \frac{2V - (\textrm{max}(V) + \textrm{min}(V))}{\textrm{max}(V) - \textrm{min}(V)}
\end{align*}
and assume $\Del$ to be uniformly distributed in the interval $[-1,1]$.

The uncertain $\vo{A}(\Del)$ and $\vo{B}(\Del)$ are expressed as:
\begin{align*}
  \vo{A}(\Del) := \sum_{i=0}^{5} \vo{A}_i \phi_i (\Del), \vo{B}(\Del) := \sum_{i=0}^{5} \vo{B}_i \phi_i (\Del)
\end{align*}
where $\phi_i(\Del)$ are $i$th-order Legendre polynomials.


The polynomial-chaos framework allows for any $\mathcal{L}_2$ function to model parametric uncertainty, not just multi-affine functions. Here, Legendre polynomials are used to capture the variation in the system matrices.

The code for simulation is publicly available\footnote{https://github.com/isrlab/LqrDiscPC}.
The components $\vo{A}_i$  and $\vo{B}_i$ can be derived from the same.

For the controller synthesis, we regulate the outputs: velocity $V$, angle of attack $\alpha$, and flight path angle $\gamma := \theta - \alpha$ defined by:
\begin{align*}
  \vo{y} := \begin{bmatrix}
    V \\
    \alpha \\
    \gamma
\end{bmatrix} = \vo{Cx}, \textrm{with } \vo{C} := \begin{bmatrix}
  1 & 0 & 0 & 0 \\
  0 & 1 & 0 & 0 \\
  0 & -1 & 1 & 0
\end{bmatrix}
\end{align*}

The optimal cost-to-go is defined using $\vo{Q}$ and $\vo{R}$ where:
\begin{gather*}
  \vo{Q} := \vo{C}^T \vo{Q}_y \vo{C}, \textrm{with } \vo{Q}_y := \textbf{diag}(10^{-1} \quad 10 \quad 10{}) \\
  \vo{R} := \textbf{diag}(10^{-4} \quad 10^{-1} )
\end{gather*}

Weights $\vo{Q}$ and $\vo{R}$ are chosen to normalize the state and control trajectories with respect to the desired peak values. The PC-based synthesis algorithm is implemented with increasing orders of approximation, up to 7.

\begin{figure}[!ht]
        \begin{centering}
   \includegraphics[width=0.45\textwidth]{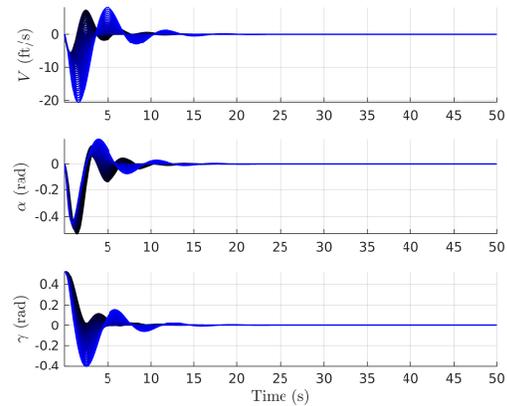}
\end{centering}
   \caption{Output trajectories obtained using the reduced order model. The initial condition is: $\vo{x}_0^T := (0\, 0\, \frac{30 \pi}{180}\, 0)$.}
   \label{yTrajRed}
\end{figure}

Figure \ref{kGain} depicts the increase in control gain with approximation order whereas figure \ref{yTrajRed} shows the closed-loop response for $V^t$, $\alpha^t$, and $\gamma^t$ with the gain $\vo{K}$ for the 7th order.
The response here is for the linear system for various values of $\Del$.
The initial condition used here is: $\vo{x}_0^T := (0\, 0\, \frac{30 \pi}{180}\, 0)$, which corresponds
to $\gamma(0) = 30^{\circ}$.
The controller is synthesized using a 7th order approximation. The colors black to blue correspond to values of $\Del$ ranging from -1 to 1.
In figure \ref{clLoopPoles}, we can see the closed-loop poles of the system created using the 7th order approximation. All poles lie within the unit circle, indicating stability of the closed-loop system. These observations are consistent with the trajectory response in figure \ref{yTrajRed}.

\section{Summary}
We derived a new algorithm for the synthesis of optimal and robust state feedback controllers for linear discrete-time systems with probabilistic parameters in their system matrices. Optimality ensures a minimum quadratic cost on the state and control action required to stabilize the system.
We followed a reduced order modeling approach that utilizes a finite-dimensional approximation built using polynomial chaos and develops a controller by solving this optimization problem exactly.
The simulations show that the controller is able to stabilize the system about the chosen operating points. 
However, it is important to note that the proposed approach works favorably only for systems proven to be stable over the entire distribution of the random parameters.
In the future, we would like to investigate whether the infinite-dimensional model in \cite{bhattacharya2019robust} could be meaningfully extended to discrete-time systems as well.
Moreover, the impact of a multi-dimensional parameter space, i.e. $\Del \subset \real^d$ with $d \geq 2$ on closed-loop stability and controller performance requires further research.

%
%
%
%
%

\bibliographystyle{elsarticle-num}
\bibliography{citations}

\vspace{100mm}

\end{document}